\documentclass[12pt]{amsart}
\usepackage{epsfig,amssymb,epic,eepic,color}
\usepackage[all]{xy}
\numberwithin{equation}{section}
 \usepackage[T1]{fontenc}  
 \usepackage{hyperref}
 
\theoremstyle{plain}    

\theoremstyle{definition}

\numberwithin{equation}{section}
\newcommand{\be}{\begin{enumerate}}
\newcommand{\ee}{\end{enumerate}}
\newcommand{\bi}{\begin{itemize}}
\newcommand{\ei}{\end{itemize}}

\def\R{\mathbb{R}}

\def\Z{\mathbb{Z}}

\def\S{\mathbb{S}}

\def\ga{\gamma}    
\def\Ga{\Gamma}

\def\al{\alpha}
\def\be{\beta}
 
\def\de{\delta}
\def\De{\Delta}
\def\vp{\varphi}

\def\si{\sigma}
\def\Si{\Sigma}

\def\ep{\varepsilon}

\def\nd{\noindent}
\def\bull{\hfill$\Box$\\}

\def\D{\mathbb{D}}

\def\p{\partial}

\usepackage{hyperref}
\usepackage{mathrsfs}
\usepackage{color,soul} 
\setstcolor{red} 
\setul{0.2ex}{0.4ex} 
\usepackage[usenames,dvipsnames]{xcolor}

\begin{document}

\vskip 1cm
\begin{center}{\sc Ren\'e Thom and an anticipated $h$-principle 
}
 \vskip 1cm
 
 {\sc Fran\c cois Laudenbach}

\end{center}
\title{}
\author{}
\address{Universit\'e de Nantes, LMJL, UMR 6629 du CNRS, 44322 Nantes, France}
\email{francois.laudenbach@univ-nantes.fr}


\subjclass[2010]{}

\thanks{This work, supported by ERC Geodycon, is an expanded version
of a lecture given in the conference: {\it Hommage \`a Ren\'e Thom}, 1-3 Sept. 2016, IRMA, Strasbourg.}

\begin{abstract}  The first part of this article intends to present the role played by 
Thom in diffusing Smale's ideas 
about immersion theory, at a time (1957) 
where some famous mathematicians were doubtful about them: {\it it is {\bf clearly}
 impossible to turn the sphere inside out!}
Around a decade later, M. Gromov transformed Smale's idea into what is now known 
as the $h$-{\it principle}. Here, the $h$ stands for {\it homotopy}.

Shortly after the astonishing discovery by Smale, Thom gave a conference 
in Lille (1959) announcing
a theorem which would deserve to be named {\it a homological $h$-principle}. 
The aim of our second part 
is to comment about this theorem which was completely ignored by the topologists in Paris, 
but not in Leningrad. We explain Thom's statement and answer the question whether it is true. The first idea
is combinatorial. A beautiful subdivision of the standard simplex emerges from Thom's article.
We connect it with the {\it jiggling} technique introduced by W. Thurston in his seminal work on foliations.

\end{abstract}

\maketitle
\thispagestyle{empty}
\vskip .5 cm

\section{From immersions viewed by Smale to Gromov's $h$-principle}
\vskip 1cm
\subsection{Thom and Smale in 1956-1957.}\label{1957}
 Important and reliable information\footnote{Curiously enough some biographies online give 1957 
 as the date of Smale's thesis despite the footnote in \cite{smale-thesis} is being quite clear on this matter.} about Smale in these years is given by M. Hirsch \cite[p.\,36]{hirsch2}:
 \begin{quote}
 I first learned of Smale's thesis at the 1956 Symposium on Algebraic Topology in Mexico City.
 I was a rather ignorant graduate student at the University of Chicago, Smale was a new PhD from
 Michigan ...
 I thought I could understand the deceptively simple geometric problem Smale addressed: {\it Classify 
 immersed curves in a Riemannian manifold.}
 \end{quote}
 Ren\'e Thom gave an invited lecture at the same Symposium. Probably, it was the first occasion 
 for Thom and Smale to meet. Let us continue reading Hirsch \cite{hirsch2}:
 \begin{quote}In the Fall of 1956, Smale was appointed Instructor at the University of Chicago.
\end{quote}

On January 2, 1957, Smale submitted an abstract to the Bulletin of the American Mathematical Society
which was published in the issue of May 1957 \cite[Abstract 380{\it t}]{smale-abstract}. This is a 14-line piece\footnote{This abstract is not included in \cite{smale-collected}.} 
titled: {\it A classification of immersions of the 2-sphere} where Smale announces\footnote{The complete article
following this announcement is \cite{smale-trans}.}
 a complete classification
of immersions the 2-sphere valued in  $C^2$ manifolds  of dimension greater than two. He wrote:
 \begin{quote}
 For example any two
$C^2$ immersions of $S^2$ in $E^3$ are regularly homotopic.
\end{quote}

\begin{figure}[h]
\includegraphics[scale=.8]{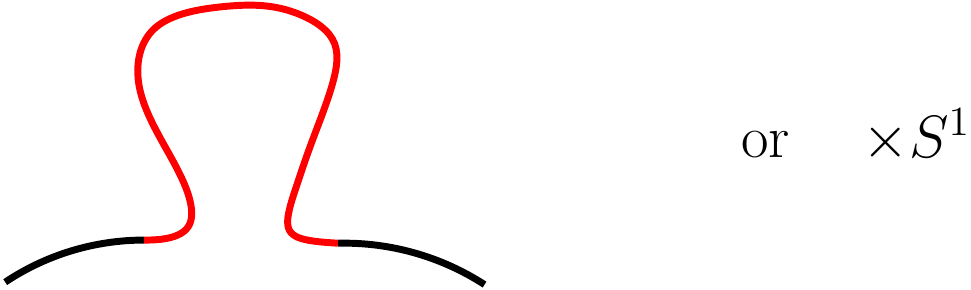}
\caption{ Corrugation in dimension one and two.}
\label{corrugation}
\end{figure}

In Spring 1957, 
Thom spent a semester as invited Professor at the University of Chicago. He spoke with Smale for hours
until he had a full understanding of Smale's ideas on immersions.  Back in 
France, Thom reported on Smale's work in a Bourbaki seminar of December 1957 \cite{thom-bourbaki} (or \cite[p.\,455-465]{oeuvres}).
It is remarkable that the written version of Thom's lecture contains the very first figure
 which had appeared in the theory of immersions. This is just a {\it bump}; according to \cite{thurston},
 W. Thurston would have called this picture a {\it corrugation} (Figure \ref{corrugation}). 
 
 I should say that the theory of corrugation is still very lively 
  for constructing concrete  $C^1$ isometric embeddings (see V. Borrelli \& {\it al.} \cite{borrelli}).
 \begin{figure}[h]
\includegraphics[scale=0.12]{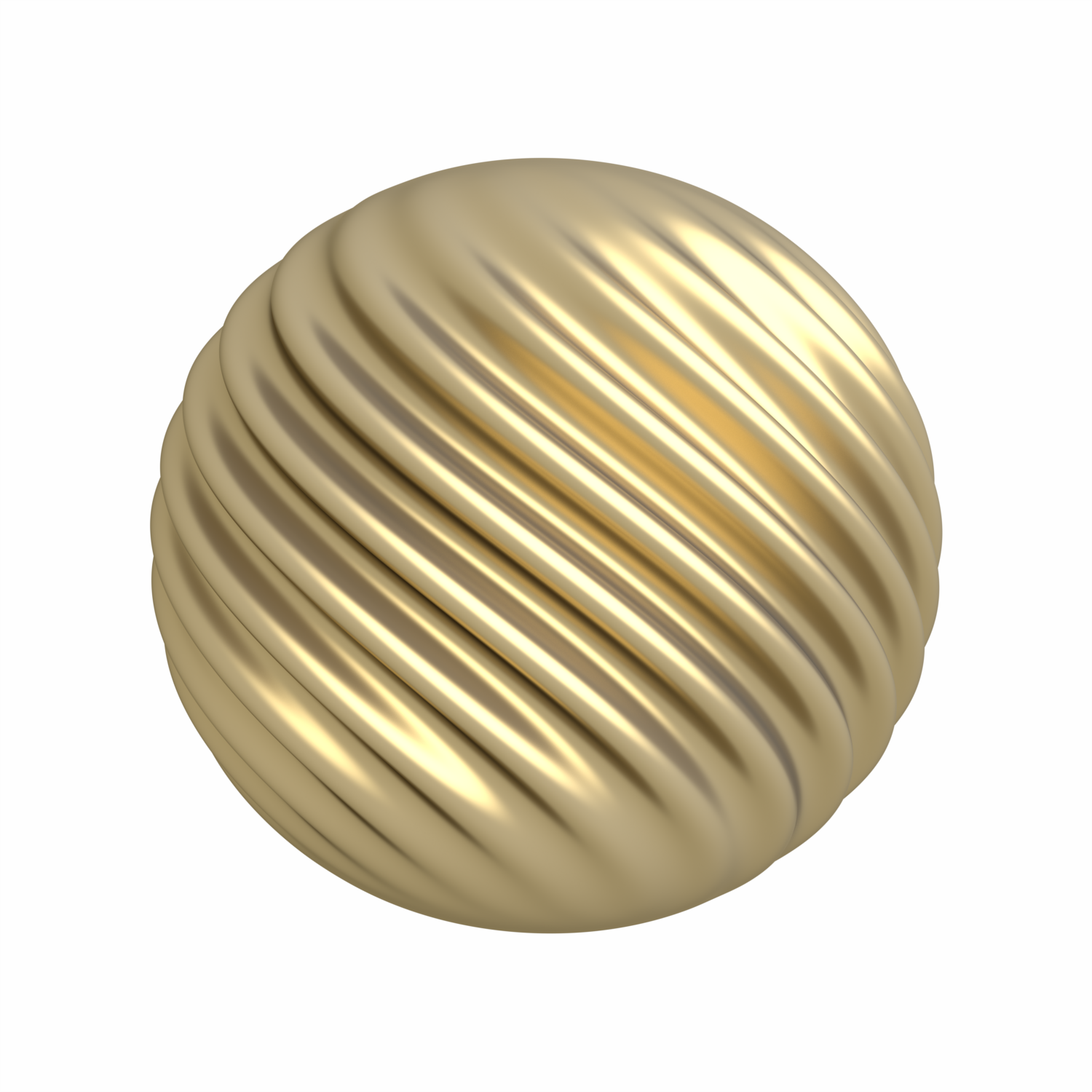}
\caption{First corrugating step for an isometric embedding of the unit sphere into the ball of radius 1/2.
By courtesy of the Hevea Project.}

 \label{corrugated-sphere}
\end{figure}

In the rest of Section 1, 
I would like to present Smale's ideas, starting from the basics, and connect them with more recent ideas.

\subsection{Immersions.}\label{immersions} 
\begin{figure}[h]
\includegraphics[scale=.7]{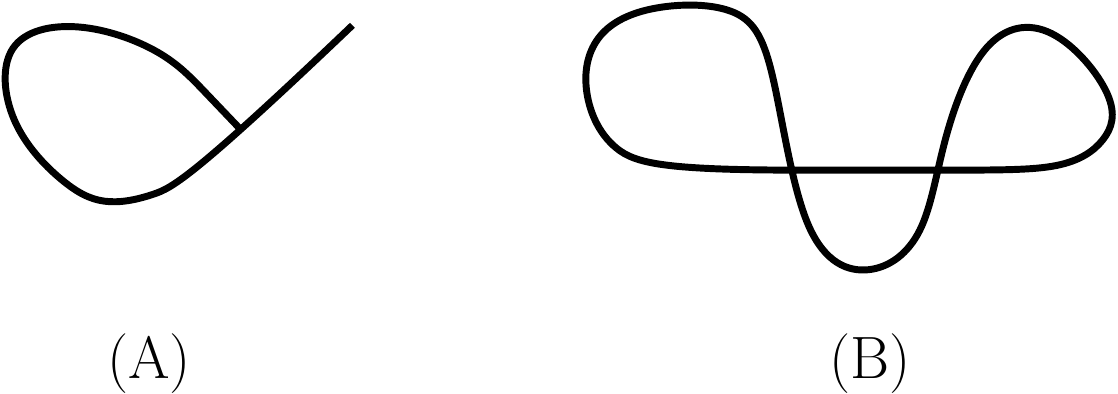}
\caption{(A) shows an embedding $(0,1) \to \R^2$ whose image is not a submanifold.
(B) shows an immersion $\S^1\to\R^2$ which does not extend to an immersion of the 2-disc.}
\label{non-embed}
\end{figure}

Given two smooth manifolds $X$ and $Y$ 
where the dimension of $X$ is not greater than the dimension of $Y$,
a $C^1$-map $f:X\to Y$ is said to be an {immersion} if its differential $df$ is of maximal rank at every 
point of $X$. An immersion can have double points but no singular points like {\it folds}. 
An immersion with no double points is said to be an {\it embedding}. In that case, the image of $f$
is a submanifold under some properness condition, more precisely when $f$ is proper (in the topological sense\footnote{Given two topological spaces $A$ and $B$,  a  continuous map from
$A$ to $B$ is said to be {\it proper} if the preimage of any compact set of $B$
is a compact set of $A$.})  from $X$ to some open subset of $Y$ (Figure \ref{non-embed}(A)).

The space of immersions from $X$ to $Y$, denoted by $Imm(X,Y)$,  is an open set in $C^1(X,Y)$
if the space of $C^1$ maps is endowed with the so-called {\it fine Whitney topology}. When $X$ is compact,
there is no concern: a sequence $\left(f_n\right)$ is convergent if and only if both sequences
 $\left(f_n(x)\right)$ and $\left(df_n(x)\right)$ converge uniformly. In what follows, we shall only consider 
 immersions whose source is compact.
 In that case, the set of immersions
 is locally contractible.
  
Two immersions $f_0,f_1: X\to Y$ are said to be {regularly homotopic} if they are joined by a path 
in $Imm(X,Y)$ or equivalently if $f_0$ and $f_1$ belong to the same path component of 
$Imm(X,Y)$.

\subsection{Whitney-Graustein Theorem} 
The immersions from the circle to the plane were classified by Whitney up to regular homotopy
in the mid-thirties \cite{whitney}. The classification reduces to the degree of the Gauss map
$$ \begin{array}{rcl}
G: \S^1&\to& \S^1\\
x&\mapsto & df_x(\p_\theta)/ \Vert df_x(\p_\theta)\Vert\,,
\end{array}
$$
where $\p_\theta$ stands for the unit tangent vector to the circle $\S^1: = \R/2\pi\Z$.
The reason why this theorem is named {\it Whitney-Graustein Theorem} is given by Whitney himself
in a footnote on page 279 of his article:
 \begin{quote}This theorem, together with a straightforward proof, was suggested to me by 
W.~C. Graustein.
\end{quote}
It is worth noticing that there is an interesting proof  of the Whitney-Graustein Theorem 
given by S. Levy in \cite[p.\,33\,-\,37]{levy} following  Thurston's idea of corrugation. 

It would be wrong to think that this classification ends the story of immersion of the circle to the plane.
A much more difficult question is the following: {\it Which immersions extend to an immersion of the disc
to the plane? For such an immersion, how many extensions are there?} An obvious necessary condition 
for a positive answer to the first question is that the degree of the Gauss map be equal to one.
 But that condition is
not sufficient as Figure \ref{non-embed}(B) shows. Actually, these questions have been 
solved by S. Blank
in his unpublished thesis. Luckily, V. Poenaru reported\footnote{It is worth noticing that Poenaru's report 
contains the drawing of 
the so-called J. Milnor's example, 
that is an immersion of the circle into the plane having two extensions to the disc 
which are not equivalent up to homeomorphism of $D^2$.} on Blank's thesis in a Bourbaki 
seminar \cite{poe2}.
The analogous questions for immersions of the $n$-sphere into $\R^{n+1}$ can be raised and remain essentially open.

\subsection{The key proposition in Smale's thesis.} Let $f_0$ denote the standard embedding of $S^2$
in $\R^3$. Choose an equator $E$ on $\S^2$, a base point $p\in E$ and two hemispheres respectively
named the northern and the southern hemisphere $H_N$ and $H_S$. We consider the space of pointed 
immersions 
$$Imm_p(\S^2,\R^3):=\left\{ f:\S^2\looparrowright \R^3 \mid f(p)=f_0(p)\text{ and }df(p)=df_0(p)\right\}.
$$
The spaces $Imm_p(H_S,\R^3)$ and $Imm_p(H_N,\R^3)$ are  defined similarly. 
The space of immersions $H_N$ to $\R^3$ whose 
1-jet $j^1\!f(x): =\bigl(f(x),df(x)\bigr) $ 
coincides with $j^1\!f_0(x)$ 
at every point $x\in E$ 
 is denoted by $Imm_E(H_N,\R^3)$. Finally, $\widetilde{Imm}_p(E,\R^3)$ stands for the space
of immersions of $E$ to $\R^3$ enriched with a 2-framing along $E$  which is fixed at $p$ and
 whose generated plane field
 is tangent to $E$ 
and standard at $p$.
The space of pointed immersions of the 2-disc to $\R^3$ is known to be contractible thanks to 
the  Alexander's  contraction which reads in the present setting: 
  $$(x,t)\mapsto p+\frac 1 t [f\bigl(p+t(x-p)\bigr)-f(p)]$$
   where $p$ lies in the boundary of $\D^2$ and $(x,t)\in \D^2\times (0,1]$. 
   When $t$ goes to 0,
 the limit of the above expression, uniformly in $x$ in the 2-disc, is the affine map $x\mapsto p+df(p)(x-p)$.\\

\nd {\sc Proposition.} ${}$

\nd 1) {\it The restriction map 
$Imm_p(\S^2,\R^3)\to Imm_p(H_S,\R^3)
$ is a Serre fibration. Its fibre over $f_0$ is homeomorphic to $Imm_E(H_N,\R^3)$.}

\nd 2) {\it The 1-jet map along the equator, $Imm_p(H_N,\R^3)\to \widetilde{Imm}_p(E,\R^3)$,
is a Serre fibration. 
Its fibre over $(j^1\!f_0)\vert_E$ is also homeomorphic to $Imm_E(H_N,\R^3)$}.\\

A map $\rho: X\to Y$ 
between two arcwise connected spaces  is said to be a  {\it Serre fibration}
when it has the parametric 
{\it Covering Homotopy Property}. More precisely, for every $\ga:[0,1]\to Y$
and every $x_0$ in $X$ with $\rho(x_0)= \ga(0)$, there exists a lift $\tilde\ga:[0,1]\to X$ of $\ga$ starting 
from $x_0$; and similarly in families with parameters in the $n$-disc. In that case, there is a long exact sequence in homotopy.

It is worth noticing that similar statements for one-dimensional source were already
present in Smale's thesis 
(published in \cite{smale-thesis}).\\

\begin{figure}[h]
\includegraphics[scale=.45]{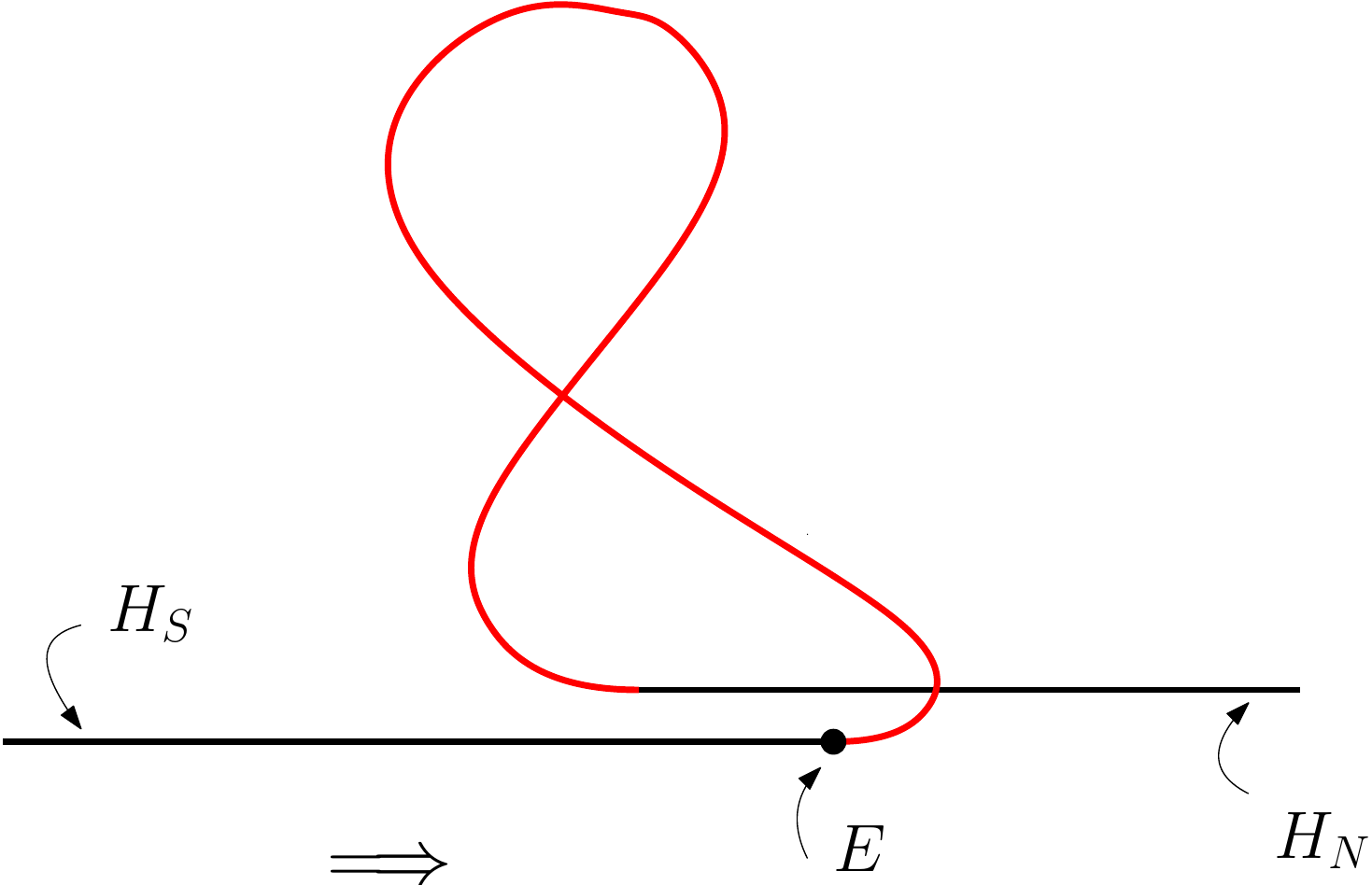}
\caption{
 The southern hemisphere moves to the right. 
Only a  northern collar of the equator in $H_N$
 is deformed.}
\label{flex}
\end{figure}

The proof of the first item is sketched by a picture which shows the {\it flexibility} that the statement 
translates (Figure \ref{flex}).
\medskip

\nd{\sc Corollary.} {\it We have $\pi_0\bigl(Imm_p(\S^2,\R^3)\bigr)=0$, that is, the space of pointed 
immersions of $\S^2\to \R^3$ is arcwise connected.}\\

\nd {\sc Proof.} Since the base of the first Serre fibration is contractible, we have 
$\pi_0\bigl(Imm_p(\S^2,\R^3)\bigr)\cong \pi_0\bigl(Imm_E(H_N,\R^3)\bigr)$.
By the second Serre fibration whose total space is contractible, we have
 $\pi_0\bigl(Imm_E(H_N,\R^3)\bigr)\cong\pi_1\bigl(\widetilde{Imm}_p(E,\R^3)\bigr)$.
 Arguing similarly for the enriched immersions of $E$ whose equator is a 0-sphere, we get 
 $$\pi_1\bigl(\widetilde{Imm}_p(E,\R^3)\bigr)\cong \pi_2\bigl(\widetilde{Imm}_p(S^0,\R^3)\bigr)
 \cong\pi_2(\{2\text{-frames in }\R^3\}\cong\pi_2(SO(3))=\pi_2(\S^3)=0.
 $$
${}$ \hfill \bull

\subsection{Concrete eversion of the sphere.} I do not intend to explain  the history of this matter.
I just give a list of references in chronological order
and add a few comments: A. Phillips \cite{phillips}, G. Francis \& B. Morin
\cite{morin}, Francis' book \cite{francis} and finally the  text and video by S. Levy \cite{levy}.

The first idea, due to Arnold Shapiro, is to pass through  {\it Boy's surface}, here noted $\Si$,
an immersion of the 
projective plane into the 3-space. Since the projective plane is non-orientable,
a tubular neighborhood $T$ of  $\Si$ is not a product. Therefore, $T$ is bounded by 
an immersed sphere $\tilde\Si$. It turns out that $\tilde\Si$ is endowed with 
the involution which consists of exchanging the two end points in each fibre of $T$. This  is realized 
by the regular homotopy 
$$ \tilde\Si\times[0,1] \to T,\ (x,t)\mapsto  x(1-2t)$$
 where the product in the right hand side is associated
 to the affine structure of the fibre of $x$. If the two faces of $\tilde\Si$ are painted with different colors,
 this move has the effect of changing the color which faces  Boy's surface. It remains to connect 
 the standard embedding of $S^2$ to Boy's surface by a regular homotopy in order to get an eversion
 of the sphere. 
 
 Remembering a walk with  Nicolaas Kuiper when he explained this construction to me, I had the feeling that he played  himself a role in it. I did not know more until very recently, when Tony Phillips 
 informed me about an article of Kuiper where his argument is written explicitly\footnote{The approaches by Shapiro and Kuiper were contemporary. As far as I know,
 nothing indicates some relationship between them.}
  \cite[p.\,88]{kuiper}.
 The video \cite{levy} does not follow the same idea: it 
 goes the way of Thurston's corrugations and is not optimal in number of multiple points of multiplicity 
 3 or more.

\subsection{Hirsch's definitive statement.} The general statement in {\it homotopy theory} of 
immersions is due to M. Hirsch \cite{hirsch}. He considers any pair $(X,Y)$ 
of  smooth manifolds.
For simplicity, assume $X$ is connected. 
The main assumption is that $\dim X \leq\dim Y$, the equality being allowed  only 
when $X$ is not closed (if $X$ is compact its  
boundary  must be non-empty).

If $f: X\to Y$ is an immersion, we have a diagram
\[ \xymatrix@R=0.3cm{
   TX\ar[dd]\ar[r]^{df}&  TY\ar[dd]\\
                                  & \\
 X\ar[r]^{f}&Y
}\]
where $df$ is a fibre bundle map {\it over} $f$ (between the total spaces of the respective tangent bundles)
which is fibrewise linear  and injective. 

Although the following terminology has been in use since Gromov's thesis only, we are going to use it here. A {\it formal immersion} is a diagram
\[ \xymatrix@R=0.3cm{
   TX\ar[dd]\ar[r]^{F}&  TY\ar[dd]\\
                                  & \\
 X\ar[r]^{f}&Y
}\]
where $f$ is only assumed to be continuous and $F$ is a fibre bundle map which is fibrewise linear 
 and injective. In the language of jet spaces, this is just a section of the 1-jet bundle $J^1(X,Y)$
 over $X$ valued in the open set of 1-jets whose linear part is of maximal rank. 
 
 With this vocabulary at hand, Hirsch's theorem states the following:\\
 
 \nd {\sc Theorem (Hirsch \cite{hirsch}).}  {\it The space $Imm(X,Y)$ of immersions from $X$ to $Y$
 has the same homotopy type\footnote{In the literature on this topic, one generally speaks of the same
 {\it weak} homotopy type, meaning that the map under consideration induces an  {\it isomorphism} 
 of  homotopy groups only (for every base point). Actually,  R. Palais 
 \cite[Theorem 15]{palais} tells us
  that the two notions are equivalent for the topological spaces we are dealing with.\label{weak}}
 as the space $Imm^{\rm formal}(X,Y)$ of formal immersions. }\\
 
 \subsection{Phillips' work on submersions.}\label{phillips}
  When the dimension of $X$ is greater than the dimension of 
 $Y$ it is natural to consider 
 submersions, that is maps of maximal rank. When such maps exist they form a space that we denote 
 by   $Subm(X,Y)$. Using again the current terminology,
 a {\it formal submersion} is  a section of the 1-jet bundle $J^1(X,Y)$
 over $X$ valued in the open set of 1-jets whose linear part is of maximal rank. Phillips'  submersion 
 theorem sounds similar to Hirsch's  immersion theorem with, nevertheless, a fundamental difference: 
 the source needs to be an open manifold. Notice that the circle has no submersion to the line 
 despite the existence of a formal submersion; a similar claim holds for any parallelizable manifold like a compact Lie
  group. \\
  
  \nd {\sc Theorem (Phillips \cite{phillips-sub}).} {\it If $X$ is an open manifold, 
  $Subm(X,Y)$ and 
$Subm^{\rm formal}(X,Y)$ have the same homotopy type.}\\

Since a foliation is locally defined by a submersion onto a local transversal, the next theorem
can be viewed as an extension of the previous one. Let $\mathcal F$ be a smooth foliation 
of the manifold $Y$. Denote its normal bundle by $\nu(\mathcal F)$; it is a vector bundle on $Y$
whose rank equals the codimension of $\mathcal F$. Denote by $\pi: TY\to\nu(\mathcal F)$ 
the linear bundle morphism over $Id_Y$ whose kernel is the sub-bundle of $TY$ made of
the tangent vectors to $Y$ which are tangent to the leaves of $\mathcal F$.  

A smooth map $f:X\to Y$ is said to be transverse to $\mathcal F$ if the bundle morphism 
$\pi\circ df:TX\to \nu(\mathcal F)$ over $f$ is fibrewise surjective. In that case, the preimage of 
$f^{-1}(\mathcal F)$ is a foliation of the same codimension as $\mathcal F$ and its normal bundle
is the pull-back $f^*\bigl(\nu(\mathcal F)\bigr)$.  We denote by $C^{\pitchfork\mathcal F}(X,Y)$
the set of smooth maps transverse to $\mathcal F$.

Given a bundle morphism $F: TX\to TY$   over $f:X\to Y$, the pair $(f,F)$ is said to be 
 formally transverse to $\mathcal F$ if $\pi\circ F$ is fibrewise surjective. By abuse, one says also that 
 $f$ is formally transverse to $\mathcal F$.\\

\nd{\sc Theorem (Phillips \cite{phillips-fol}).} {\it The space $C^{\pitchfork\mathcal F}(X,Y)$
has the same homotopy type as the space of maps which are formally transverse to $\mathcal F$.}\\

\nd {\sc Remark.} All previous theorems reduce the understanding of immersions, submersions or maps transverse to foliations
 from the homotopic point of view to the understanding of the 
corresponding formal problems. And the latter reduces to classical homotopy theory: the matter is to find 
 sections to some maps and thus it reduces to {\it well-known obstructions}. This does not mean that the 
 homotopy type of the formal spaces in question is computable. In general it is not, as the homotopy 
 groups of the spheres are  not  completely computable.
 
 The aim of Gromov's approach which we are going to describe below is to consider all previous problems 
 as particular cases of a general {\it principle}. \\

 \subsection{Differential relations after M. Gromov.} \label{open-h} The main reference here is Gromov's book
 \cite{gromov}. A simplified approach is described in Y. Eliashberg \& N. Mishachev's book
 \cite{eliash-misha}; the new tool is their {\it holonomic approximation Theorem} which was first proved 
 in \cite{holonomic}. 
 
The preface of  \cite{eliash-misha} starts as follows:
\begin{quote}
 A {\it partial differential relation} $\mathcal R$ is any condition imposed on the partial derivatives 
 of an unknown function.
 \end{quote}
 If the unknown function in question is a smooth map from $X$ to $Y$ -- we limit ourselves to this 
 case\footnote{More generally, the unknown could be a section of a given smooth bundle over $X$.} --
  a simple definition  consists of saying  that $\mathcal R$ is a subset in a jet space\footnote{Since
   this is going to be forgotten, I recall that the concept of jet space is due to Charles Ehresmann.}
   $J^r(X,Y)$ for some integer $r$. Recall that in coordinates an element of this jet space is just the data 
   of a point $a\in X$,
  a point $b\in Y$ and a Taylor expansion of order $r$ at $a$ with constant term $b$.
  
  The expression {\it $h$-principle} comes from the article of Gromov \& Eliashberg \cite{removal}; they write:
  \begin{quote}The principle of weak homotopy equivalence for {\it etc.}
  \end{quote}
  Later on, this expression is abbreviated to {\it h-principle}\footnote{I already 
  commented on the word {\it weak} in 
   footnote \ref{weak}.  Concerning the word {\it principle}, I feel uncomfortable 
   with a principle which is not always true, and worse, whose domain of validity remains unknown.
   That means $h$-principle is not a gift from heaven.}. With these authors we say that the 
   {\it parametric $h$-principle holds for $\mathcal R$} if the inclusion 
   $$sol\, \mathcal R\hookrightarrow sec\, \mathcal R
   $$
  is a homotopy equivalence between 
  $$sol\, \mathcal R:= \{f:X\to Y\mid j^r\!f \text{ is valued in } \mathcal R\}
 $$
  and 
   $$sec\,\mathcal R:=\{\text{sections of } J^r(X,Y)\text{ valued in }\mathcal R\}\,.$$

 One can think of an element of $\sec\mathcal R$ as a {\it formal} solution of the problem posed by
 $\mathcal R$. A section valued in $\mathcal R$ which is of the form $j^r\!f$ is said to be {\it holonomic}
  or {\it integrable}. The integrablity is prescribed by the vanishing on the section in question of a list 
  of 1-forms (called a {\it Pfaff system})
   which are naturally defined on the manifold $J^r(X,Y)$. For instance, when $r=1$ and $Y= \R$,
  a section $s$ is integrable if and only if its image is {\it Legendrian} for the canonical {\it contact form} $\al$
  which reads $dz- \sum p_idq^i$ in canonical coordinates, that is, if $s^*\al=0$.
 \medskip

 \nd{\sc Theorem (Gromov \cite{gromov-thesis}).} {\it 
 If $\mathcal R$ is an open set in $J^r(X,Y)$ which is invariant under the natural right action of 
 $\text{Diff}(X)$ and if $X$ is an {\sc open} manifold (meaning that no connected component is closed),
 then the parametric $h$-principle holds true for $\mathcal R$.}\\
 
 The proof also goes through corrugations as said for Smale's theorem. Of course, the corrugations
 are not developed in the range; there is no room for corrugating. They are developed in the domain.
 This is very clearly explained in Eliashberg-Mishachev's book \cite{eliash-misha}.\\

 \nd{\sc Remark.} Another very important condition on a differential relation leads to an $h$-principle; it is 
 when the relation is {\it ample}. In that case $X$ does not need to be an open manifold. Here,  the 
 $h$-principle follows from the famous {\it convex integration} technique which was invented by Gromov in
  \cite{convex} (see Gromov's book \cite{gromov}).  A complete account on this is given in D. Spring's
 book \cite{spring}. The end of  Eliashberg-Mishachev's book \cite{eliash-misha} focuses on  convex 
 integration applied to the $C^1$ isometric embedding problem (Nash-Kuiper);
   Borrelli \& {\it al.} \cite{borrelli} converted their theoretical result into an algorithm richly illustrated by 
   pictures of $C^1$ {\it fractal objects}, as the authors say. Despite the great interest of the subject,
   I do not intend to enter more deeply into it 
as it is less connected to the work of Thom than what follows.\\
 
  Going back to the, say {\it open}, $h$-principle stated above, one sees that the  previously mentioned 
   results by Smale, Hirsch and Phillips are clearly covered by Gromov's theorem. One could be 
   disappointed that only the
  1-jet space is involved. The simplest way to find new examples with differential relations of higher order
 consists of the following construction, which naturally appears in Thom's singularity theory \cite{thom-sing}
 as it is shown in the next subsection.
 
 Consider a proper submanifold $\Si\subset J^{r-1}(X,Y)$ or a proper stratified set with nice 
 singularities (for instance, with {\it conical singularities} in the sense of \cite{app-bismut}); the important point 
 is that transversality to any stratum $\Si_i\subset \Si$ implies transversality to all other strata in some 
 neighbourhood of $\Si_i$ in $ J^{r-1}(X,Y)$. Assume that $\Si$ is {\it natural}, that is invariant under the action 
 of ${Dif\!f(X)}$. The transversality to $\Si$ is obviously a differential relation of order $r$.
 This differential relation which we denote by $\mathcal R_\Si$ is open
 and invariant under the action of ${Dif\!f(X)}$. Thus, if $X$ is open, Gromov's theorem applies.
 
 \subsection{\bf Examples from singularity theory.}\label{ex} 
 
 For a first concrete example, 
 take $\dim X=\dim Y=2$ and consider the stratified set $\Si$ of 1-jets of rank less 
 that  2. It is made of two strata: one stratum is the set of jets of rank 1; it has codimension 1
 and is denoted by $\Si^1$ in the so-called Thom-Boardman notation \cite{boardman}. The other stratum  is the 
 rank zero one; it has codimension 4 in our setting. Their union is a stratified set 
 with conical singularities which is natural and proper. Thus $\mathcal R_\Si$ satisfies the open
 $h$-principle if $X$ is open.\\
 
  \begin{figure}[h]
\includegraphics[scale=.45]{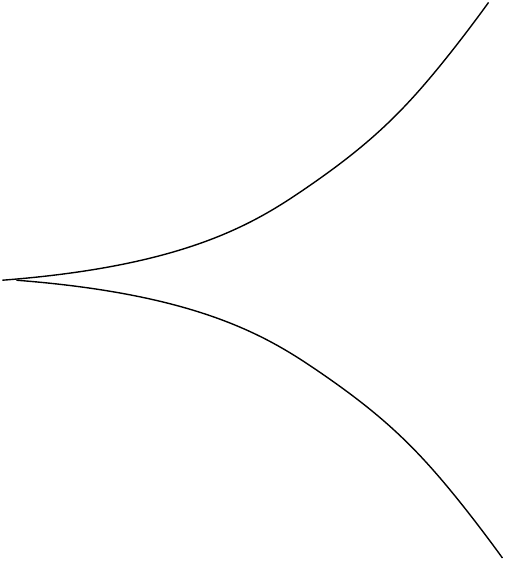}
\caption{ Local image of a cusp in a 
 two-dimensional manifold.}\label{cusp}
\end{figure}
 The next example leads to an order 3 differential relation. One starts with the first example 
 and looks at a 2-jet $\xi\in \mathcal R_\Si$; say it is based in $a\in X$. 
 Since $\xi$ is transverse to $\Si$, it does not project to the zero 1-jet. Therefore, 
 it determines
 the tangent space in $a$ to the {\it fold locus} $L\subset X$ where the rank of any germ of map $f$
 realizing $\xi$ is exactly 1.  In our setting, $L$ is one-dimensional. On the other hand, $\xi$ determines
 the kernel $K_a$ of the differential $df_a$. 
 Thus, there is a natural 
 stratification of $\mathcal R_\Si\subset J^2(X,Y)$: one stratum is $\Si^{1,0}$ which is made of 2-jets 
 where $K_a$ is transverse to $T_aL$; the second one, denoted by $\Si^{1,1}$, is made of 2-jets 
 where $K_a$ is tangent to $T_aL$. The stratum $\Si^{1,0}$ is an open set in $\mathcal R_\Si$
 and $\Si^{1,1}$ has codimension 2 in $J^2(X,Y)$; it is a conical singularity of $\mathcal R_\Si$.
 Thus, if a 3-jet  is transverse to $\mathcal R_\Si$, it is the jet of a germ having an isolated {\it cusp}
 from which emerge two branches of fold locus (see Figure \ref{cusp}).\\
 
 \subsection{\bf Thom's transversality theorem in jet spaces.}
 This was exactly the subject of Thom's lecture at the 1956 Symposium in Mexico City
  that I mentioned at the very beginning of this piece.  
  Incidently, this  theorem  will play a fundamental role in singularity theory, as the above discussion 
  lets us foresee. The statement is the following:\\
  
  \nd{\sc Theorem (Thom \cite{mexico}).} {\it Let   $\Si$ be a submanifold in a $r$-jet bundle
  $E^{(r)}\to X$
  over a manifold $X$. Then, generically\footnote{A property is said to be generic in a given topological space $F$ (here, it is the space of integrable sections with the $C^0$ topology or the Whitney topology 
  evoked in Subsection \ref{immersions}) if it is satisfied by all elements in a residual subset (that is, 
  an intersection of countably many open dense subsets).}, 
  an integrable section  of $E^{(r)}$ is transverse to $\Si$.}\\
  
  This theorem is remarkable in two ways:
  
  \nd 1) The usual transversality statement tells us that any section of $E^{(r)}$ can be approximated
  by a section transverse to $\Si$. But, the integrability condition is a {\it closed 
  constraint}\footnote{The space of integrable sections is closed with empty interior in the space of all 
  sections.} and even if we started with an integrable section, the transverse approximation could be 
  non-integrable.

  \nd 2) The same proof, by inserting the given map
  in a large family of maps which is transverse to $\Si$ as a whole, 
  works both for the usual transversality theorem and for the transversality theorem  with constraints.
  \medskip
  
  For many years I tried to understand whether the statements of Thom and Gromov 
   were somehow related. For instance, does the $h$-principle hold for the relations $\mathcal R_\Si$ 
   from subsection \ref{ex}? 
   The answer was shown  to be {\sc no}  in general,
     in a note with
  Alain Chenciner \cite{chencine}. 
  Quoting from its abstract: 
  \begin{quote}
   A section in the 2-jet space of Morse functions is not always homotopic to a holonomic section.
\end{quote}
${}$\\

\section{Integrability and related questions}

 \subsection{\bf Thom's point of view in 1959.}
 The title of the lecture given by Ren\'e Thom at the  1959 conference organized by the CNRS 
 in Lille (France) is striking when compared with the terminology that would appear ten years later:
 \begin{quote}
 Remarques sur les probl\`emes comportant des in\'equations diff\'erentielles globales
 \end{quote}
 which I translate into:
 {\it Remarks about problems involving global differential inequations.}\\

\nd The setting is the same as in Gromov's theorem from Subsection \ref{open-h}
and, for consistency with what precedes,  $\mathcal R$ still denotes 
an open set in the jet space $J^r(X,Y)$, except that now the openness of $X$ is not assumed. There are 
two chain complexes naturally associated with $\mathcal R$:
\begin{enumerate} 
\item $C_*(\mathcal R)$ is the complex of continuous\footnote{Replacing  {\it continuous} with {\it smooth}
changes the complex to a {\it quasi-isomorphic} subcomplex, meaning that the homology is unchanged.}
singular simplices.

\item $C_*^{\rm int}(\mathcal R)$ is the subcomplex generated by the differentiable simplices valued in 
$\mathcal R$ which are integrable (or holonomic) in the sense
that each 1-form from the integrability Pfaff system
vanishes on them.
\end{enumerate}

Here, a $k$-simplex is  a map from the standard $k$-simplex $\De^k\subset \R^{k+1}$ to $\mathcal R$.
Thanks to the so-called {\it small simplex} Lemma, up to quasi-isomorphism, it is sufficient to consider 
holonomic smooth simplices of the form: $\si= j^r f\circ \underline\si$ where $\underline\si$ is a $k$-simplex 
of the base $X$ and $f$ is a smooth map defined near the image of $\underline \si$ with values in $Y$.\\

\nd {\sc Theorem (Thom \cite{thom-ineq}.)} {\it 
The inclusion $C_*^{\rm int}(\mathcal R)\hookrightarrow C_*(\mathcal R)$ induces 
an isomorphism in homology for $*<\dim X$ and an epimorphism for $*=\dim X$.}\\

For instance, if $X$ is closed and $s: X\to \mathcal R$ is a section, then the cycle $s(X)$ 
(at least with $\Z/2\Z$ coefficients when $X$ is non-orientable) is homologous to a holonomic 
{\it zig-zag}, that is a cycle of the form $j^rf(X)$ where $f:X\to Y$ is multivalued.\\

\subsection{\bf What happened afterwards.} This article was actually  only
an announcement. The proof of 
the theorem was outlined in three pages, and was difficult to read 
although  some ideas were visibly emerging;
for instance the {\it sawtooth}, which is an antecedent to the {\it jiggling} intensively used by Thurston
in the early seventies \cite{thurston-fol}. No complete proof ever appeared.
Unfortunately, the report by Smale 
in the Math. Reviews \cite{smaleMR}
was somewhat discouraging for anyone who would have tried
 to complete Thom's proof. Here is the final comment of this report:
 \begin{quote}
 \{The author has said to the reviewer that, although he believes his proof to be valid 
 for $r= 1$, there seem to be further difficulties in case $r>1$.\}
 \end{quote}
 
 Nevertheless, David Spring has known for a few years that Thom's statement holds true (see his note
 \cite{spring-note}).
  His unpublished proof  is based on the holonomic approximation theorem of Eliashberg \& Mishachev
 \cite{holonomic} when $*<\dim X$. In the remaining case, he also needs Poenaru's foldings
  theorem \cite{folding}. I should say that the holonomic approximation theorem 
  is in germ in Thom's announcement; his horizontal sawtooth is closely related to the construction made in \cite{holonomic}.
  
  When reading Thom's article for preparing the edition\footnote{The team of editors of Thom's works
  was initiated by Andr\'e Haefliger and is directed by Marc Chaperon.}
   of his collected mathematical works 
  \cite{oeuvres}, I was no more able 
  to complete the proof in the way indicated by Thom,
  but I discovered a beautiful  object 
  in that article.
  I  first translate the original few lines  into English 
    and then, in the next subsection, I shall state the lemma which I could extract from these lines.
  
 \begin{quote}
 [The proof] mainly relies on the construction of a deformation (homotopy operator) from 
 the complex of all singular differentiable simplices to the integrable simplices. Such a deformation 
 has to be \guillemotleft \ hereditary \guillemotright, that is, compatible with the restriction to faces. 
Moreover, as the problem is local in nature, it will be sufficient  to construct this deformation for an open set
in $J'^r(\R^n,\R^p)$.

... ...
 
 Let $b^k$ be a  $k$-dimensional simplex, $b^n$ an $n$-dimensional simplex,
  $n\geq k$; let $b'^k$
 be a subdivision of $b^k$ and $s$ a simplicial map from this subdivision $b'^k$ of $b^k$
 to $b^n$. The finer the subdivision $b'$  is, 
 the more  the map $s$ has a \guillemotleft\,strong gradient\,\guillemotright\  in the sense that the quotient $[s(x)-s(y)]/(x-y)$, for every pair of points
 $x,y\in B^k$ close enough, becomes larger and larger.
 \end{quote}
  
Here, the question is: why do such a subdivision and simplicial map exist?\\

\subsection{\bf Thom's subdivision.}
\begin{figure}[h]
 \begin{center}
\includegraphics[scale=.6]{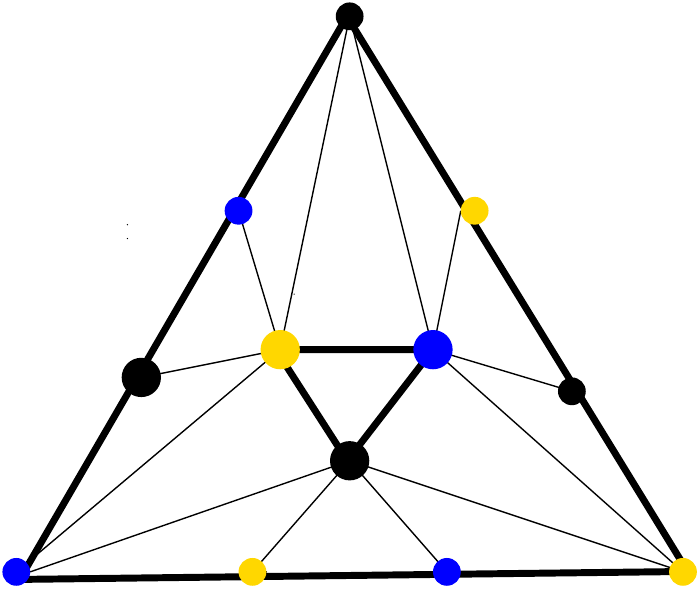}
\end{center}
\caption{Thom subdivision in dimension two; $s_2$ is defined by the coloring. }
\label{subdivision2}
\end{figure}

Here is the statement that I cooked up for translating the preceding lines 
into a more precise language\footnote{I was recently informed by Michal Adamaszek that this subdivision is known to theoricists of computing and is now called {\it the standard chromatic subdivision}. A similar figure to Figure \ref{subdivision2} appears in the article by M. Herlihy and N. Shavit \cite{herlihy}.}.\\

\nd {\sc  Lemma.} {\it There exists a sequence $\left( K_n,s_n\right)_n$, where $K_n$ 
is a linear subdivision of $\De^n$
and $s_n: K_n\to \De^n$ is a {\bf simplicial} map such that:
\begin{enumerate}
\item ({\bf Non-degeneracy}) for each  $n$-simplex $\de_n\subset K_n$, the restriction $s_n\vert_{\de_n}$ is surjective; 
\item ({\bf Heredity}) for each $(n-1)$-face $F$ of $\De^n$ we have:
\begin{center} 
$$
\left\{
\begin{array}{l}
F\cap K_n\ \cong\  K_{n-1}\,,\\
s\vert_F\ 
\cong\  s_{n-1}\,.
\end{array}
\right.$$
\end{center}
\end{enumerate}
}
 Here, the symbol $\cong$ stands for simplicial isomorphism; if a numbering of the vertices of $\De^n$
 is given there is a canonical simplicial  isomorphism $F\cong \De^{n-1}$ for every facet $F$. The non-degeneracy 
 somehow translates Thom's strong gradient condition.  
 \begin{figure}[h]
 \begin{center}
\includegraphics[scale=.6]{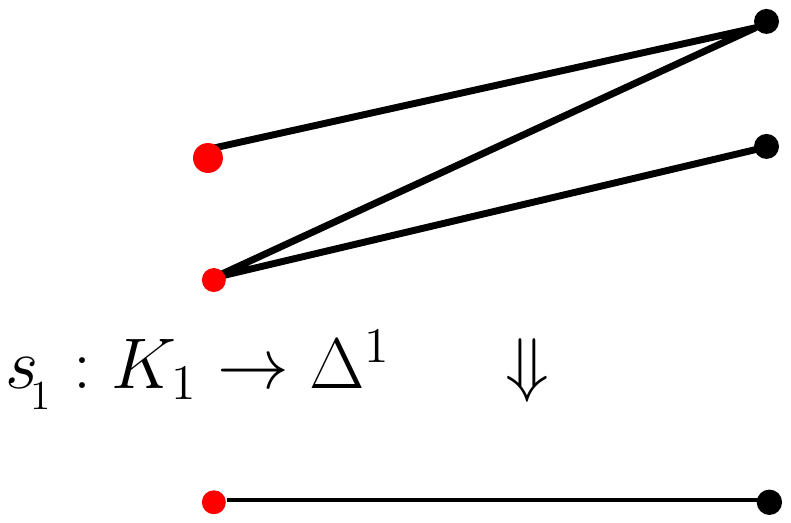}
\end{center}
\caption{Thom subdivision in dimension 1 and its folding map.}
\label{subdivision1}
\end{figure}
The proof can be obtained by induction on $n$ in the way which is  illustrated 
by passing from Figure \ref{subdivision1} to
Figure~\ref{subdivision2}: put a small $n$-simplex $\de^n$ 
{\it upside down} in the interior of $\De^n$
and join each vertex of $\de^n$ to the facet of $\De^n$ lying in front of it which is already 
subdivided by induction hypothesis.

One can think of $s_n$ as a folding map from $\De^n$ onto itself. Due to the heredity property, we have:
\begin{itemize}
\item Any polyhedron can be folded onto itself.
\item The folding can be iterated $r$ times: 
$$
\begin{array}{lcl}
K_n^{(r)}&=& \left(s_n^{(r-1)}\right)^{-1}(K_n)\\
s_n^{(r)}&= &s_n\circ s_n ^{(r-1)}
\end{array}
$$
\end{itemize}
Notice that the folding map of any order is endowed with an hereditary {\it unfolding} 
homotopy to Identity. \\

\subsection{\bf Jiggling formula} It is now easy to derive  a natural {\it jiggling} formula, 
 using the same terminology as Thurston's in \cite{thurston-fol}, but without any measure consideration.

Equip $X$ with a Riemannian metric. Let $DX\to X$ be a tangent disc bundle such that the 
{\it exponential} map $\text{exp}: DX\to X$
is a submersion. Choose a triangulation $T$ of $X$ finer than the open covering 
$\left\{ \text{exp}_x(D_xX)\mid x\in X\right\}$. Fix an integer $r$. The $r$-th jiggling map 
is the section of the tangent bundle defined by
$$
\begin{array}{c}
j^{(r)}: X\to DX,\\
j^{(r)}(x)= \text{exp}_x^{-1}\bigl(s_n^{(r)}(x)\bigr)\,.
\end{array}
$$
This map is piecewise smooth. Moreover, 
the larger $r$ is,  the more {\bf vertical}   the jiggling is. 
\begin{figure}[h]
 \begin{center}
\includegraphics[scale=.6]{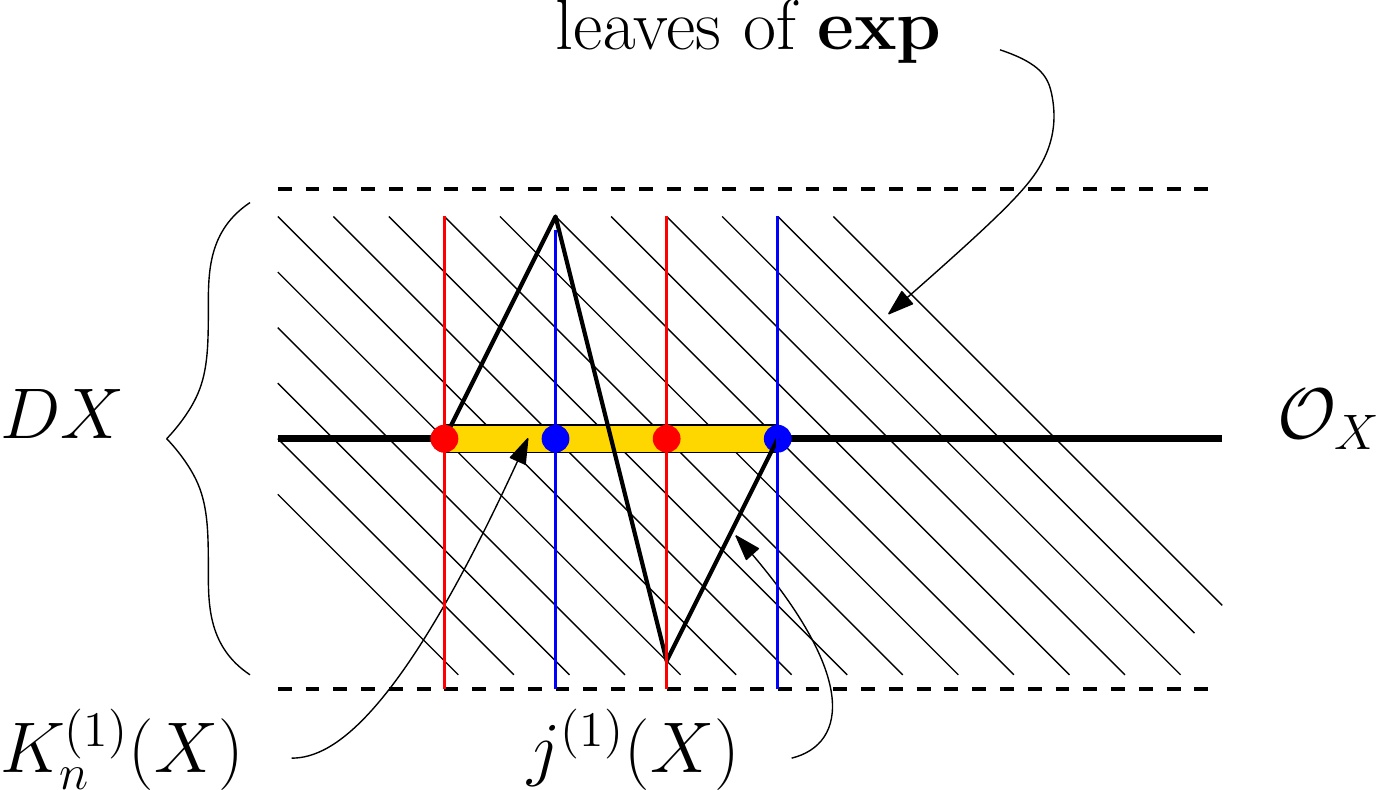}
\end{center}
\caption{The jiggling map of order $r=1$. The vertical lines are fibres of $TX$.}
\end{figure}
As a consequence, for $r$ large enough, $j^{(r)}(X)$ is {\it quasi-transverse} to the tangent space to the 
exponential foliation $\mathcal F_{\rm exp}$, that is, for any simplex $\tau$ of the $r$-th 
Thom subdivision
of  $T$ the smooth image $j^{(r)}(\tau)$ shares no tangent vector with the tangent space to the leaves 
of $\mathcal F_{\rm exp}$. Actually, when $X$ is compact, 
this quasi-transversality  holds 
with respect to any compact family of $n$-plane fields which are 
 transverse to the fibres of $TX$, in place of 
$\mathcal F_{\rm exp}$.

\subsection{\bf Going back to immersions.} \label{going-back}
This is part of a joint work with Ga\"el Meigniez \cite{l-m}.

First, recall that one can reduce oneself to consider only immersions of codimension 0. Indeed, any formal 
immersion $(f,F)$ from $X$ to $Y$ ($\dim X<\dim Y$) has a normal bundle; it is the vector bundle 
over $X$ which is the cokernel $\nu(f,F)$ of the monomorphism $TX\to f^*TY$ through which 
$F$ factorizes. Thus, immersing $X$ to $Y$ is equivalent to immersing a disc bundle of $\nu(f,F)$ to $Y$
and the latter is a codimension 0 immersion.

In what follows, we assume that  $X$ is compact with non-empty boundary and has the same dimension
 as $Y$. For free, a formal immersion $(f,F)$ from $X$ to $Y$ gives rise to a foliation 
 $\mathcal F_X:= F^{-1}(\mathcal F_{{\rm exp}_Y})$ which foliates a neighbourhood of the zero section 
 $O_X$
 of $TX$. Indeed, since $F$ maps fibres to fibres surjectively, $F$
 is transverse to the exponential foliation of $Y$ (defined near the zero section 
 $O_Y$ of $TY$). 
 
  \begin{figure}[h]
 \begin{center}
\includegraphics[scale=.6]{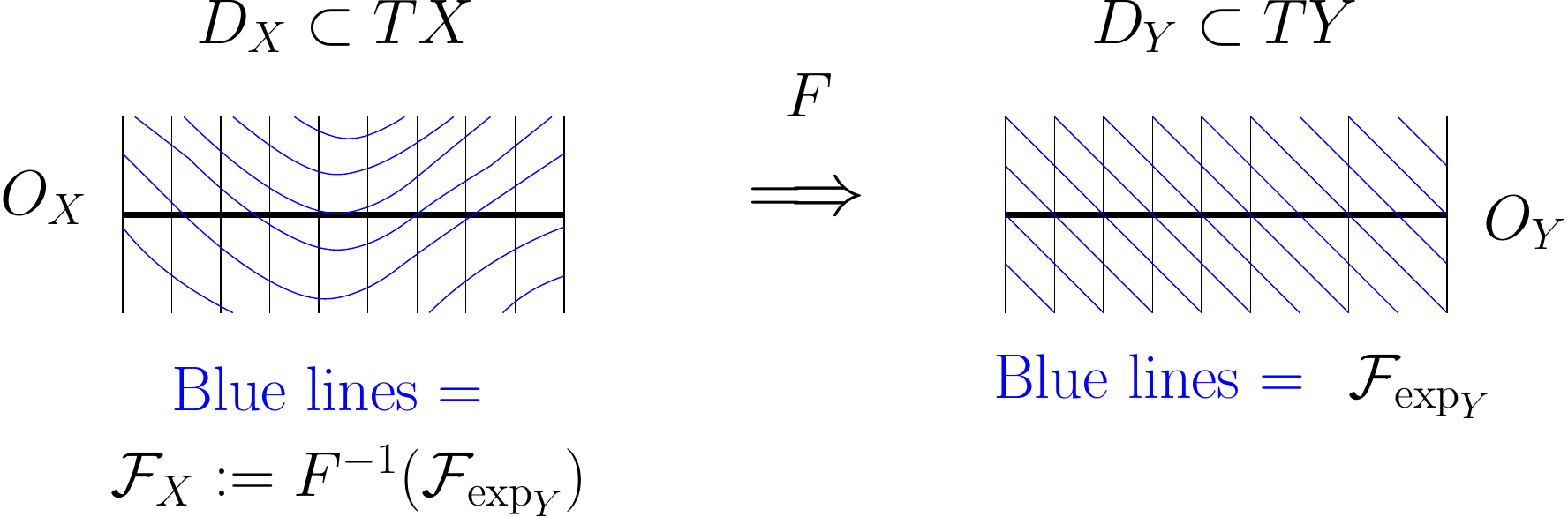}
\end{center}
\caption{The left part shows a tangential Haefliger structure.}
\label{formal}
\end{figure}

 Such a (germ of)
 foliation like $\mathcal F_X$ is called a {\it tangential Haefliger structure} 
 or a {\it $\Ga_n$-structure} on $X$.
 We refer to \cite{haefliger} for more details on this important notion.
Since there is no reason for $\mathcal F_X$ to be transverse to  $O_X$, the trace of 
 $\mathcal F_X$ on $X= O_X$ is in general a {\it singular} foliation. 
 
 Actually, those singularities 
 are responsible for the flexibility associated with that concept: they allow for operations like  
 induction (or pull-back) and homotopy (or  concordance).
  Let us emphasize that a $\Ga_n$-structure is mainly 
 a $\check {\rm C}$ech cocycle of degree one valued in the groupoid of germs of diffeomorphisms of 
 $\R^n$. This allows one to induce such a structure on a polyhedron or a CW-complex. 
 A concordance between two $\Ga_n$-structures $\xi_0,\xi_1$ on $X$ is just a $\Ga_n$-structure on 
 $X\times[0,1]$ which induces $\xi_i$ on $X\times\{i\},\,i=0,1$. There is a classifying space $B\Ga_n$
 in the following sense:
  the $\Ga_n$-structures on $X$, up to concordance, are in 1-to-1 correspondance
   with the homotopy classes $[X,B\Ga_n] $, as for vector bundles.
  
  In our setting, the Haefliger structure in question is enriched with a {\it transverse geometric structure}
  invariant under
   holonomy: each 
  transversal to $\mathcal F_X$ is endowed with a submersion to $Y$  which
  is preserved when moving the transversal by isotopy along the leaves (this point being obvious since the 
  leaves in question are contained in the inverse images of points in $Y$);
   such a $\Ga_n$-structure will be named a 
  $\Ga_n^Y$-structure. 
  In particular, if $O_X$ were transverse 
  to $\mathcal F_X$,
  then $X$ would be endowed with a submersion to $Y$, that is  an immersion to $Y$ as $\dim X= \dim Y$.
  Therefore, the aim is to remove the singularities of the $\Ga_n^Y$-structure, that is, to find a 
  {\it regularizing 
  concordance} of $\Ga_n^Y$-structures from $\mathcal F_X$ to a $\Ga_n^Y$-structure whose underlying 
  foliation is transverse to the zero section. 
  
  In the next subsection we give a brief review of the regularization  problem, and in the last subsection
  a sketch of the regularization is given in our setting of immersions in codimension 0 of compact manifolds 
  with non-empty boundary and no closed connected components.
  
  \subsection{\bf About the regularization of $\Ga$-structures.} Let $\xi$ be a $\Ga_q$-structure on 
   an $n$-dimensional manifold $X$. In general, the underlying foliation $\mathcal F(\xi)$ is supported
  in a neighbourhood of the zero-section in a vector bundle $\nu(\xi)$
  of rank $q$, called the normal bundle to $\xi$.
 This normal bundle remains unchanged along a concordance. If $\xi $ is regular, that is, if 
 $\mathcal F(\xi)$ is transverse to the 0-section of $\nu(\xi)$,
 then the trace of $\mathcal F(\xi)$ on $X$ is a genuine foliation 
 whose normal bundle is canonically isomorphic to $\nu(\xi)$.
 Therefore, a necessary condition to be regularizable is that $\nu(\xi)$ embed into the tangent 
 bundle
 $TX$\footnote{By abuse, we confuse a vector bundle with its total space.}; in particular, $q\leq n$.
 
 Andr\'e Haefliger was the first to prove that any $\Ga_q$-structure on an {\bf open} manifold $X$
 whose normal bundle embeds into $TX$ is regularizable \cite{haefliger70} (or \cite[p.\,148]{haefliger}).  
 That follows from two things:  first, the classifying property of the classifying space $B\Ga_q$: the latter 
 is equipped with a universal $\Ga_q$-structure which induces by pull-back all others; 
 second,  the Phillips transversality theorem to  a foliation
 \cite{phillips-fol} (see the statement in Subsection \ref{phillips}). Today, this regularization theorem 
 is frequently referred to as the Gromov-Haefliger-Phillips theorem.
 
 The next step was done by W. Thurston \cite{thurston-fol}. If $q>1$, even when $X$ is closed,
 any $\Ga_q$-structure satisfying the necessary condition is regularizable. The case $q=n$
 is the toy case. The only technique is the famous {\it jiggling} lemma whose proof is 
 quite tricky in terms of measure theory, despite its author considered it as obvious. 
 Exactly at this point, our jiggling based on the Thom subdivision is much simpler; 
 moreover, it works in families.
 
 The final step is the codimension-one case for closed manifolds, a piece of work indeed. 
 Generally it is known in the following form:\\
 
\nd{\sc  Theorem (Thurston \cite{thurston-one}).} {\it Every hyperplane field is homotopic to the field tangent to some codimension-one foliation}.\\
 
  Actually, the main part of that result is a regularization theorem for $\Ga_1$-structures. 
  In addition to the jiggling technique,
 there are many subtle points (simplicity of the group of diffeomorphisms, intricate constructions, {\it etc.}).

 \subsection{\bf  Regularization of transversely geometric $\Ga_n$-structures.}
\begin{figure}[h]
 \includegraphics[scale= 0.8]{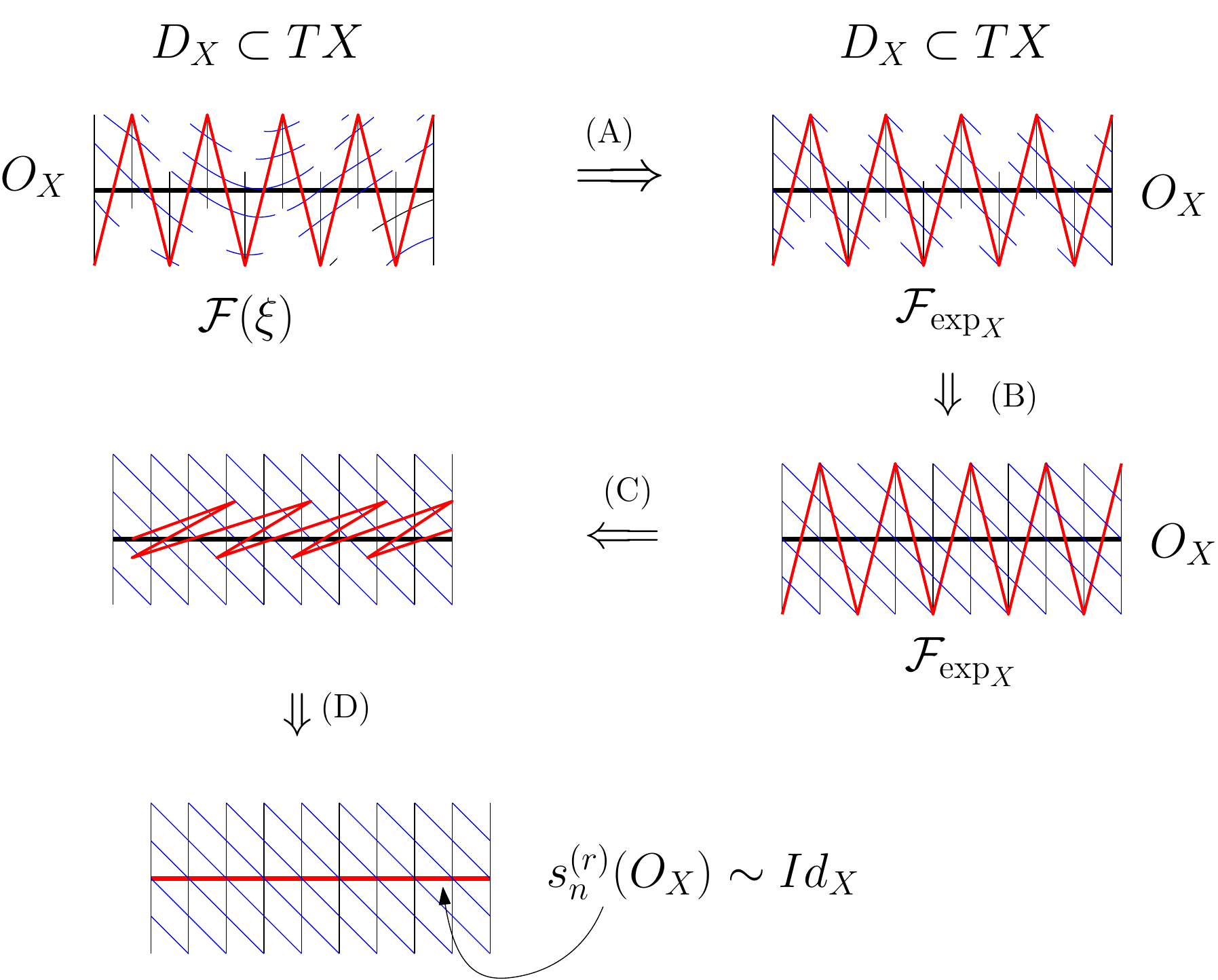}
 \caption{The scheme of the regularization in four steps.}
 \label{formal2}
 \end{figure}

 In Subsection \ref{going-back}
 we reduced the problem of immersion to a problem of regularization of some 
 $\Ga_n^Y$-structure $\xi$ on $X$ associated with the given formal immersion and shown in 
 Figure \ref{formal}. The exponent $Y$ reminds us that we are considering a $\Ga_n$-structure 
 endowed with some transverse geometry which here consists of being endowed with 
 a submersion to $Y$. The scheme shown in Figure
 \ref{formal2}, and  on which I am going to comment, summarizes an ordinary regularization
  (which would work even if $X$ were closed).
  It will appear in the end  that this regularization is easily enriched  with a transverse geometric structure
  when $X$ is open. It is worth noticing that the problem is the same whatever the transverse geometry is.
 In place of {\it submersion to $Y$} one could have a symplectic or contact structure, a complex structure
 or a codimension-one foliated structure {\it etc.}. For any geometry\footnote{We take the concept 
 of geometry in the sense of Veblen \& Whitehead \cite{veblen_31}
 which could be rewritten in the more modern language of sheaves.}, the regularization is the same.
 
 First, the jiggling is chosen, meaning that the order of the Thom subdivision $r$ is fixed once and for all. 
 This $r$
 is chosen so that $j^{(r)}(X)$ is quasi-transverse to the following codimension-$n$ foliations or plane fields:
 \begin{itemize}
 \item  the foliation $\mathcal F(\xi)$ underlying the given $\Ga_n$-structure $\xi$ (this foliation was denoted by $\mathcal F_X $ in the particular case of Figure \ref{formal});
 \item the exponential foliation $\mathcal F_{\text{exp}_X}$; 
 \item every $n$-plane field which is a  barycentric combination\footnote{Recall that the
 space of $n$-planes tangent to the total space $TX$ at $(x,u)$, $x\in X,\, u\in T_xX$, 
 and transverse to the {\it vertical tangent space} (that is, the kernel of $D_{(x,u)}\pi$ where $\pi: TX\to X$ denotes the projection) 
 is an affine space.} of the two previous ones.
 \end{itemize} 
 The homotopy from the zero-section to $j^{(r)}(X)$ gives rise to an obvious concordance which is not
 mentioned in the scheme of Figure \ref{formal2}.
 
 Step (A) is exactly Thurston's concordance in \cite{thurston-fol}. By using the above-mentioned 
barycentric combination, some generic $(n+1)$-plane field $\Pi$
is chosen on $TX\times[0,1]$ quasi-transverse to $j^{(r)}(X)\times[0,1]$. Since the trace of $\Pi$
on each simplex of the jiggling is 0- or 1-dimensional, such a trace is integrable. Thus, a $C^0$ 
approximation of $\Pi$ is integrable in a neighbourhood of $j^{(r)}(X)\times[0,1]$. This gives the 
concordance (A) and explains the reason why some part of the tube $D_X$ has been deleted
 from the initially foliated domain.
 
 Step (B) is just the inclusion using the fact that the exponential foliation exists on the whole tube.
 Step (C) uses the interpolation $\text{exp}^t$, $t\in[0,1]$, from $Id_{DX}$ to $\text{exp}: DX\to X$
 given by:
 $$(x,u)\longmapsto \left(\text{exp}_x(tu), \left(D_{(x,tu)}\text{exp}_x(tu)\right)(1-t)u\right).$$
  It allows one to slide $j^{(r)}(X)$ along the leaves of the exponential foliation keeping the 
  quasi-transversality to each simplex\footnote{For the reader who does not like complicated formulas, 
  I suggest a more topological approach of the previous interpolation. Let $U$ be a {\it nice} tubular 
  neighbourhood of the diagonal $\De$ in $X\times X$; here, nice means that the two projections 
  $\pi^v$ and $\pi^h$ respectively given by $(x,y)\mapsto (x,x)$ and $(x,y)\mapsto (y,y)$ are $n$-disc 
  bundle maps. If $U$ is small enough, a Riemannian metric provides an identification of $\pi^v$ with a
  tangent disc bundle to $X$. In that case, $\pi^h$ is the corresponding exponential map. Hence, 
  the mentioned 
  interpolation is just a contraction of the fibres of $\pi^h$.}. 
  Observe that the vertical homothety does not have such a property.
  When $t= 1$, we finish with the folding map $s_n^{(r)}(X)\to X$.
  Step (D) is just the unfolding of $s_n^{(r)}$, that is, its hereditary homotopy to $Id_X$. Again, at each time 
  of the homotopy, the image polyhedron (contained in the zero-section)
   is quasi-transverse to the exponential foliation. This finishes the regularization of $\xi$ as a 
   $\Ga_n$-structure. In general, it is not  possible to extend the transverse geometry to the concordance. 
   But,
   this is possible when $X$ is an open manifold as we are going to explain\footnote{Only the idea of the 
   proof is given here. For more details we refer to \cite{l-m}.}.\\
   
   If $X$ is an $n$-dimensional manifold without closed connected component, endowed with a triangulation
   $T$, there exists a {\it spine}, that is, a subcomplex $K$ of dimension $n-1$ such that, for any 
   neighbourhood $N(K)$, there is an isotopy of embeddings $\vp_t:X\to X$ whose time-one 
   maps $X$ into $N(K)$ (see for instance \cite[p.\,40-41]{eliash-misha}). 
   
   Restricting ourselves 
   to $K\subset X$, let us consider the concordance  $(W,\mathcal F_W)$ of $\Ga_n$-structures
   obtained by 
   concatenation and time
   reparametrization of the four concordances 
    described right above from $j^{(r)}(K)$ to $K\subset O_X$. Here, $W\subset TX\times [0,1]$
    is  piecewise linear homeomorphic to $X\times [0,1]$; and $\mathcal F_W$
    is a  codimension-$n$ foliation  defined near $W$
    and transverse to the fibres of $ TX\times [0,1]\to  X\times [0,1]$ which induces $\mathcal F(\xi)$
    over $t=0$ and $\mathcal F_{{\rm exp}_X}$ over $t=1$. Moreover, $\mathcal F_W$ is 
    quasi-transverse to every simplex of $W$. Therefore, since $W$ is $n$-dimensional, every leaf 
    meets each simplex of $W$
    in one point at most\footnote{Here, it is necessary to make a jiggling in the time direction.
    The cell decomposition of $W$ is then prismatic (simplex$\times$interval). 
    Each prismatic cell has a Whitney triangulation (canonical up to the numbering of the vertices of $X$) 
    \cite[Appendix II]{whitney-prism}.}.
        
    By construction, 
    $W$ {\it collapses} onto its initial face $W_0:= j^{(r)}(K)$.
    We recall that  a simplicial complex $W$
   collapses to $W_0$ if there is a sequence of elementary collapses $W_{q+1}\searrow W_{q}$
   starting with $W$ and ending with $W_0$. An {\it elementary collapse}
    means that $W_{q+1}$ is the union of 
   $W_{q}$ and a simplex $\si_q$ so that $\si_q\cap W_{q}$ consists of the boundary of $\si_q$ with an
   open facet  removed. 
   The elementary collapse $W_{q+1}\searrow W_{q}$
   gives rise to an {\it elementary isotopy} $\chi_q^t$ pushing $W_{q+1}$ into itself,
   keeping $W_q$ fixed, and ending with $\chi_q^1(W_{q+1}) $ as close to 
   $W_q$ as we want. Due to the quasi-transversality to $\mathcal F_W$ this isotopy extends 
   to a neighbourhood of $\si_q$ as a {\it foliated isotopy} 
   $\tilde \chi_t$, meaning that leaf is mapped to leaf 
   at each time.
   \begin{figure}[h]
   \includegraphics[scale=0.8]{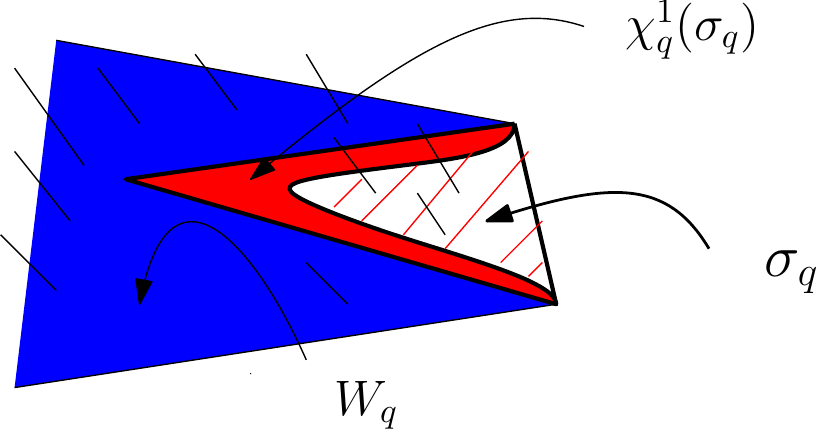}
   \caption{A few leaves of $\mathcal F_W$ are drawn in black.
    The simplex $\si_q$ is the union of the red full part 
  and the red dashed part. }
   \label{collapse}
  \end{figure}
   
   By induction on $q$, assume that the transverse geometric structure already exists on
   the foliation $\mathcal F_W\vert_{W_q}$. Then, by pulling back through  $\tilde\chi_q^1$,
   this structure extends to the foliation  $\mathcal F_W\vert_{W_{q+1}}$. Finally, the whole 
   foliation $\mathcal F_W$ is enriched with the considered geometry, for instance a submersion to $Y$.
   And hence, $N(K)$ is endowed with a submersion to $Y$.\\

\subsection {\bf Sphere eversion again.} The main advantage of this proof 
based on the Thom subdivision and its associated jiggling
is that it works in families (or with parameters). It is sufficient to choose the order $r$
large enough so that a common jiggling is convenient for each member of the family.

For instance, if $f_0$ denotes the inclusion $S^2\hookrightarrow \R^3$ and $f_1:= -f_0$, these two
 immersions are formally homotopic\footnote{I learnt this very simple formula from Ga\"el Meigniez.} by:
$$(f_t,F_t): \ (x,\vec u)\mapsto \left(tf_1(x)+(1-t)f_0(x),R_{Ox}^{\pi t}(\vec u)\right).
$$
 Here, $t\in[0,1]$ is the parameter of the homotopy, $x$ is a point in $\S^2$ and $\vec u$ is a vector in 
 $\mathop{\R}\limits^{\to}\,\!\!^3$, the vector space 
 underlying the affine space $\R^3$,
 tangent to $S^2$ at $x$; and $R_{Ox}^{\pi t}$
 stands for the Euclidean rotation of angle $\pi t$ in $\mathop{\R}\limits^{\to}\,\!\!^3$
   around the oriented axis directed by $\vec x$. When $t= 1$, we have indeed $F_1(x,-)= d_xf_1(-)$, 
 the differential of $f_1$ at $x$.
 
 By thickening, we have a one-parameter family $F_t$
 of formal submersions of $S^2\times (-\ep,+\ep)$ to $\R^3$. Thus, we have a one-parameter family
 of $\Ga_3$-structures equipped with a transverse geometry (the local submersion to $\R^3$). 
 The regularization by the Thom jiggling method -- one jiggling for all foliations
 $F_t^{-1}(\mathcal F{{\rm exp}_{\R^3}}$) -- gives rise to a one-parameter family of submersions
 $S^2\times (-\ep,+\ep) \to \R^3$ joining the respective  thickenings 
 of $f_0$ and $f_1$. The restriction to $S^2\times\{0\}$ is a regular homotopy from $f_0$ to $f_1$.
  This is the desired sphere eversion.
 
 Here is a final remark. Since $f_0$ and $f_1$ have the same image,
 we get that the space of {\sc non-oriented}
 immersed 2-spheres in the 3-space is not simply connected. Maybe, those who were skeptical 
 about the sphere eversion thought that the orientation should be preserved. Of course, if an orientation
 is chosen on the initial sphere, it propagates along any regular homotopy. But, as the image is changing
 it does not prevent us  from a change of orientation when the final image is the same as the original one.
 This is a   phenomenon of {\it monodromy} well-known for detecting  non-simply-connectedness.

\vskip 1cm

\nd{\bf Acknowledgements.}  I am indebted to Tony Phillips who provided important information
to me  about the time of the sphere eversion. I thank Paolo Ghiggini for helpful advice.
I also thank Peter Landweber who read very carefully the first posted version  and sent good remarks to me. \\

\end{document}